\documentclass[11pt,reqno]{amsart}
\usepackage{amsmath, latexsym, amsfonts, amssymb,
amsthm, amscd,epsfig,enumerate}
\advance\hoffset -.75cm

\usepackage{amsfonts}
\usepackage{latexsym}
\usepackage{amsmath}
\usepackage{amssymb}

\usepackage{color}

\newcommand{\Z}{\mathbb Z}
\newcommand{\N}{\mathbb N}
\newcommand{\K}{\mathbb K}

\newcommand{\E}{\mathbb E}

\newcommand{\Zd}{{{\mathbb Z}^d}}

\newcommand{\eps}{\varepsilon}

\newcommand{\cE}{\mathcal{E}}

\newcommand{\pr}{\mathbb P}

\newcommand{\ident}{{\mathchoice {\rm 1\mskip-4mu l} {\rm 1\mskip-4mu l}
{\rm 1\mskip-4.5mu l} {\rm 1\mskip-5mu l}}}

\newtheorem{theorem}{Theorem}[section]
\newtheorem{lem}[theorem]{Lemma}
\newtheorem{cor}[theorem]{Corollary}
\newtheorem{rem}[theorem]{Remark}
\newtheorem{pro}[theorem]{Proposition}

\begin{document}

\title
{A self-regulating and patch subdivided population}

\author[L.~Belhadji]{Lamia Belhadji}
\address{L.~Belhadji,  Laboratoire de Math\'ematiques Rapha\"el Salem, UMR 6085, CNRS -
Universit\'e de Rouen,
Avenue de l'Universit\'e, BP.~12, 76801 Saint Etienne du Rouvray, France.
}
\email{lamia.belhadji\@@univ-rouen.fr}
\author[D.~Bertacchi]{Daniela Bertacchi}
\address{D.~Bertacchi,  Universit\`a di Milano--Bicocca
Dipartimento di Matematica e Applicazioni,
Via Cozzi 53, 20125 Milano, Italy.
}
\email{daniela.bertacchi\@@unimib.it}

\author[F.~Zucca]{Fabio Zucca}
\address{F.~Zucca, Dipartimento di Matematica,
Politecnico di Milano,
Piazza Leonardo da Vinci 32, 20133 Milano, Italy.}
\email{fabio.zucca\@@polimi.it}

\begin{abstract}
We consider an interacting particle system on a graph which, from a macroscopic point of view,
looks like $\Z^d$ and, at a microscopic level, is a complete graph of degree $N$ (called a patch).
There are two birth rates: an inter-patch one $\lambda$ and an intra-patch one $\phi$. Once a
site is occupied, there is no breeding from outside the patch and the probability $c(i)$ of success
of an intra-patch breeding decreases with the size $i$ of the population in the site. We prove the existence
of a critical value $\lambda_{cr}(\phi, c, N)$ and a critical value $\phi_{cr}(\lambda, c, N)$.
We consider a sequence of processes generated by the families of control functions $\{c_n\}_{n \in \N}$ and degrees
$\{N_n\}_{n \in \N}$; we prove, under mild assumptions, the existence of a critical value
$n_{cr}(\lambda,\phi,c)$. Roughly speaking we show that, in the limit, these processes behave
as the branching random walk on $\Z^d$ with inter-neighbor birth rate $\lambda$ and on-site birth rate $\phi$.
Some examples of models that can be seen as particular cases are given.
\end{abstract}

\maketitle


\noindent {\bf Keywords}: contact process, restrained branching random walk, epidemic model, phase transition, critical parameters.

\noindent {\bf AMS subject classification}: 60K35.

\baselineskip .6 cm


\section{Introduction}

Stochastic models for the demographic expansion (or contraction) of a biological population
have attracted the attention of many researchers since Galton and Watson introduced the
branching process as a very simple and non-spatial model to study the survival of surnames.
After the appearance of interacting particle systems and of the paper of Harris \cite{cf:Harris74}
on the contact process it has become increasingly clear that the spatial structure of
the environment and the interaction between individuals bring a desirable complexity in the models
one can define (see for instance \cite{cf:DurrettLevin94} for a discussion on spatial models).

We want to introduce and study a model for a biological population where individuals breed and die
living on a patchy habitat.
Clearly, if we substitute the words \textit{infect} and \textit{recover} for \textit{breed} and \textit{die},
our model also serves as a picture for the spreading of an infectious disease.

\begin{figure}\label{fig:hab}
\centering
\epsfig{figure=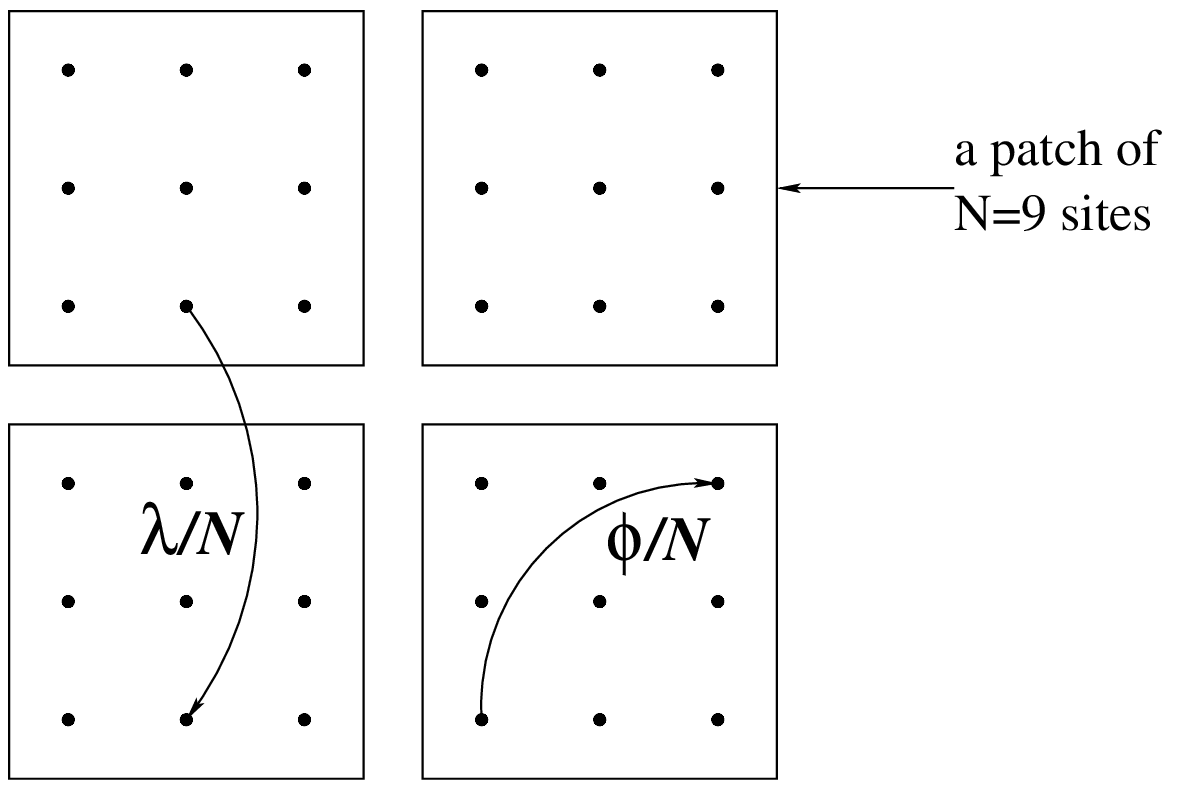,height=5cm}
\caption{The habitat.}
\end{figure}

The main idea is that the environment is subdivided in patches, the patches
are centered at the vertices of
the $d$-dimensional lattice $\Z^d$, each patch has $N$ sites
and each site may host a colony of one or more particles (see Figure \ref{fig:hab},
where each square is a patch, $N=9$ and $d=2$).
More precisely,
the patches of the habitat are complete graphs $\K_N$ with $N$ vertices
(that is any vertex is a neighbor of any other vertex, including itself). In
each single patch the sites are enumerated from 1 to $N$ and in each site there might
be a priori any number of particles.
Formally, the environment is $\cE^d_N:=\Z^d \times \K_N$ and the state space
is $\N^{\cE^d_N}$. The \textit{patch} at $x \in \Z^d$ is the set $\{(x,r): r \in \K_N\}$.
We write $(x,r) \sim (y,r_1)$, meaning that the two sites are in neighboring patches, if and only
if $x\sim y$, that is, 
$|x-y|_{\Z^d} = 1$
(where $|\cdot|_{\Z^d}$ denotes the usual distance in the graph $\Z^d$).

The dynamics of the process is the following.
Each particle dies at rate 1; when a site is vacant it receives births from each particle in the
same patch at  rate $\phi/N$ and from each particle in the neighboring patches at  rate $\lambda/N$.
When a site $x$ is occupied by $i\ge1$ particles, then it allows birth attempts
only from particles in the same patch, but an attempt succeeds only with probability
$c(i)$ (the function $c:\N\to [0,1]$ represents a self-regulating mechanism as in \cite{cf:BPZ07}).
 We suppose that $c(0)=1$ (when a site is vacant the first reproduction trial is
almost surely successful) and that $c$ is nonincreasing  
(the more individuals are present,
the less likely it is to produce new offspring).
Clearly, if $c(i)=0$ for all $i> i_0$ and $c(i_0)>0$, then in each site there
might be at most $i_0+1$ particles (see Section \ref{sec:model} for the formal
definition of the transition rates).

%
Note that we have two parameters regulating the breeding
of a particle: an intra-patch breeding rate $\phi$ and an
inter-patch breeding rate being $\lambda$.
This kind of
biparametric interaction is similar to the one introduced in
\cite{cf:Schin2002} and later modified in \cite{cf:BL}. In both papers
the population is divided in equally sized clusters, one cluster per vertex of
$\Z^d$; in \cite{cf:Schin2002} individuals in the same clusters die at the same time,
while in \cite{cf:BL} they die one at the time. In their model the cluster has no inner structure,
that is, the only important thing is whether it is empty or not (if not, immigration is not allowed).
In our model instead of clusters, we have patches and the process depends not only on
the total number of particles in the patch, but also on their position among the sites within
the patch. For instance let there be $n$ particles in the patch. If they are on the
same site, then all sites in the patch but this accept immigration (that is,
reproductions from neighboring patches) while intra-patch reproductions
are, a priori, always allowed. If they are scattered in $n$ different sites
all these sites will block immigrants.

To strengthen the dependence of the process on the number of individuals
per site, we admit a self-regulating mechanism, namely
there exists a control function $c \in [0,1]^\N$ such that if a site
is already occupied by $j$ individuals, newborns are accepted only with probability
$c(j)$. This kind of regulating mechanism has been introduced in
\cite{cf:BPZ07} as a general model of which the branching random walk
and the contact process are particular cases (when $c\equiv 1$ and $c=\delta_0$
respectively).
Note that the presence of this function $c$ mimics the fact that when
an ecosystem is already exploited by many individuals, the younglings
may not have enough food to live and to reach sexual maturity, thus
many of them will not
contribute to the breeding process. On the
other hand many animals feed on their younglings if food is scarce.


We observe that the idea of a patchy environment has already been introduced
in \cite{cf:BL2}. The main difference is that in \cite{cf:BL2} the inter-patch
interaction is allowed only between the centers of neighboring patches,
while in our model each site may receive immigrants from the neighboring patches.

In the present work we mainly study the effect of the parameters on the probability
of survival of the entire population. In particular we proved if $\lambda$ and $\phi$
are both small then there is almost sure extinction (Theorem~\ref{th:brwcoupling}).
This result is in agreement with the conclusion of the mean field model in a
particular case (see Subsection~\ref{sec:meanfield}). If $\phi$ is small
then survival is guaranteed only by $\lambda$ sufficiently large;
while $\phi$ large and a ``mild'' control function $c$ implies survival
even when $\lambda=0$ (the precise statements are Theorem~\ref{th:lambdagen} and
Corollary~\ref{cor:lambdagen} where we give explicit conditions on $\phi$ and
$c$). Analogously a large $\lambda$ guarantees survival even when $\phi=0$,
while for small $\lambda$ a large $\phi$ (and a suitable $c$) is needed to
ensure survival (see Theorem~\ref{th:phigen}).

To investigate the impact of different environment structures on
the process we let the patch dimension go to infinity. We obtain
Theorem~\ref{th:cgen}(1) which agrees with the intuition that the
limiting process is a branching random walk with intra-site breeding
rate $\phi$ and inter-site rate $\lambda$. The reader may want to compare
this convergence to the one studied in \cite{cf:BDS89},
where a long range contact process is proved to converge to a branching random walk.
We also study the effect of a progressive removing of the control function:
we obtain survival when $\phi>1$ and $c$ is close to $1$ (Theorem~\ref{th:cgen}(2)).
Finally we study the asymptotic behavior  of the critical parameter
of the general process and of some particular cases when
the size of
the patch goes to infinity (Proposition~\ref{th:Ninfty} and Corollary~\ref{th:irpk})
or when the control function is progressively removed (Proposition~\ref{th:cinfty}).

Here is a brief outline of the paper.
In Section \ref{sec:model} we introduce our process and discuss
some models that can be seen as particular cases.
Subsection~\ref{sec:meanfield} is devoted to the study of the mean field model.
In Section~\ref{sec:results} we discuss our main results; the proofs can be found in
Section~\ref{sec:proofs}.
Section~\ref{sec:open} is devoted to the discussion of open questions.

\section{The model}\label{sec:model}

Given two nonnegative parameters $\phi$ and $\lambda$,
$N\in\N$ and a nonincreasing function $c$ such that $c(0)=1$,
 the transition rates at $(x,r)$  for the process $\eta_t$ are
\begin{equation}\label{eq:genrates}
\begin{split}
0\to1 & \text{\ \ at rate\ \ }\frac{\lambda}{N}\sum_{(z,r_1):(z,r_1)\sim (x,r)}\eta(z,r_1)
+\frac{\phi}{N}\sum_{r_1 \in \K_N}\eta(x,r_1);\\
i\to i+1& \text{\ \ at rate\ \ }\frac{\phi}{N}c(\eta(x,r))\sum_{r_1 \in \K_N}\eta(x,r_1) ,\text{ for }1\le i;\\
i\to i-1 & \text{\ \ at rate\ \ }i.
\end{split}
\end{equation}
An equivalent description of the breeding mechanism is that each
particle breeds at rate $\phi$ inside its patch and at rate $\lambda$
towards each neighboring patch. Newborns choose the target sites
inside the destination patch at random. The reproduction
from a different patch is successful only if the target site is vacant,
while inside the patch the reproduction is regulated by the function $c$.

Since we are interested in the survival (with positive probability)
or extinction of the population, we always assume that the initial configuration
has a finite number of particles. In fact we observe that the process starting from
infinitely many particles gets never extinct at any time, while the process starting
from one initial particle survives with positive probability if and only if
so does the process starting from a finite number of particles. Indeed,
if it survives starting from one particle clearly it survives starting
from finitely many particles. Conversely, consider a finite initial condition
$\eta_0$ and a family of PRPs $\{\widetilde \eta_{t,i}\}_{i=1, \ldots, N}$,
where $N= \sum_{x \in \cE^d_n} \eta_0(x)$, each one starting from a single
particle (the positions are chosen according to the starting vertices
of the particles of $\eta_0$). It is not difficult to prove that
the process $\sum_{i =1}^N \widetilde \eta_{t,i}$ dominates the
original PRP (starting from $\eta_0$). Hence, if the original PRP survives
then $\sum_{i =1}^N \widetilde \eta_{t,i}$ survives, that is, $\widetilde \eta_{t,i}$
survives for some $i$. Thus there is a positive probability of survival starting
with a single particle.

In the sequel we refer to the process with these rates as the \textit{Patchy
Restrained Process} or briefly PRP. When we need to stress the dependence on the
parameters, we write PRP($\lambda,\phi,c,N$).
We observe that many explicit models are recovered as particular cases of this
patchy habitat model.

In particular, by choosing $\phi=0$ we obtain a process
which does not depend on $c$ and allows at most one particle per site.
We call this process the \textit{contact process} on $\cE^d_N$, briefly CP. It is not difficult
to prove that this process has a critical parameter
$\lambda_{CP}(\cE^d_N)
\in \left[ \frac{1}{2d},\lambda_{CP}(\Z^d)\right]
$
(namely the population dies out for $\lambda$ below it and survives with positive
probability for $\lambda$ above).
Indeed $\lambda_{CP}(\cE^d_N)\le\lambda_{CP}(\cE^d_1)=\lambda_{CP}(\Z^d)$ since
$\lambda_{CP}(\cE^d_N)$
is nonincreasing with respect to $N$ and
it is possible to prove
that $\lambda_{CP}(\cE^d_N)\downarrow 1/2d$ as $N\to\infty$ (see Proposition~\ref{th:Ninfty}(1)).

On the other hand, if we choose $c=\delta_0$, we obtain a process that we call
\textit{biparametric contact process} (or BCP($\lambda,\phi$)) on $\cE^d_N$ which coincides with the contact
process when $\phi=0$ (this is in general a process which allows at most one particle
per site and has two reproduction rates -- the inter-patch $\lambda$ and the intra-patch
$\phi$).

The \textit{individual recovery process} (briefly IRP)
introduced in \cite{cf:BL} can be seen as a particular case of the PRP, namely
by setting $N=1$ and $c=\ident_{\{0,\ldots,\kappa-1\}}$ (the patch has only one site, which
can host at most $\kappa$ particles).

\smallskip

It is natural to study sequences of PRPs (see Subsection~\ref{sec:meanfield}
and Section~\ref{sec:results})
by considering a sequence $\{N_n\}_{n\in\N}$
and a corresponding sequence of controlling functions $\mathbf{c}=\{c_n\}_{n\in\N}$
(we keep $\lambda$ and $\phi$ fixed). We are interested in the behavior
of the process as $n \to \infty$.
We assume that $c_n\le c_{n+1}$ (the regulation gets weaker as $i$ grows).
Note that $c_n(i)$ may be
positive for all $i$, thus there could be
no \textit{a priori} bounds on the site carrying capacity.
We denote by $c_\infty(i):=\lim_{n\to\infty}c_n(i)$ for all $i$.
\smallskip

We conclude this section by listing some particular cases of the PRP
which are modifications of the IRP.

\subsection{Logistic IRP}

If we choose $N=1$ and $c(i)=\max\{0,1-i/\kappa\}$
(for some $\kappa \in \N \setminus \{0\}$)
then we get a cluster of size $\kappa$ at each site in $\Z^d$
and an interaction between sites which is allowed
only when the target 
is empty.
This kind of interaction was introduced by Schinazi
\cite{cf:Schin2002} with a cluster recovery clearing
mechanism (or mass extinction), and modified
by Belhadji and Lanchier with individual recoveries in \cite{cf:BL}.
Note nevertheless that this model differs from the IRP
in the fact that the breeding inside the cluster becomes
increasingly difficult as we approach the full
carrying capacity $\kappa$.
For this process $\eta_t$, that we call the \textit{logistic} IRP,
the transition rates at site $x$ are:
\begin{equation}\label{eq:logisticIRP}
\begin{split}
0\to 1& \text{\ \ at rate\ \ }\lambda\sum_{z\sim x}\eta(z);\\
i\to i+1& \text{\ \ at rate\ \ }i\phi\left(1-\frac i\kappa\right),\text{ for }1\le i\le \kappa-1;\\
i\to i-1 & \text{\ \ at rate\ \ }i,\text{ for }i\ge1.
\end{split}
\end{equation}
Clearly the presence of the logistic factor $1-i/\kappa$ in the breeding
rate may be interpreted as a self-regulating mechanism
of the population which slows down the reproductions
when the patch is almost completely exploited.

\subsection{Self-regulating IRP}

If $N=1$ in general we get a process $\eta_t$ that we call \textit{self-regulating} IRP.
The process has the following rates at site $x$:
\begin{equation}\label{eq:genB}
\begin{split}
0\to1 & \text{ at rate }\lambda\sum_{z\sim x}\eta(z);\\
i\to i+1& \text{ at rate }i\phi c(i);\\
i\to i-1 & \text{ at rate }i.
\end{split}
\end{equation}

\subsection{Logistic IRP with persistent inter-patches reproduction}\label{sec:logistic}

Given a general PRP, we may identify all the sites in the same patch,
meaning that we consider the process
$\xi_t(x)=\sum_{r\in\K_N}\eta_t(x,r)$. We refer to $\xi_t$ as the
projection on $\Zd$ of $\eta_t$.
The projection is a Markov process only when $c\equiv\delta_0$.
In this case it has the following  transition rates at site $x$:
\[
\begin{split}
i\to i+1& \text{\ \ at rate\ \ }\left(i\phi+\lambda\sum_{z\sim x}\xi(z)\right)\left(1-\frac iN\right),\text{ for }0\le i\le N-1;\\
i\to i-1 & \text{\ \ at rate\ \ }i,\text{ for }i\ge1.
\end{split}
\]
We note that this process is similar to the IRP of \cite{cf:BL} with cluster size
equal to $N$. Nevertheless it has two main differences: the interaction between clusters 
is always active (while in the IRP it vanishes once the cluster is
nonempty) and it gets more difficult to increase the number of individuals
in the cluster if the cluster is crowded.

\subsection{The mean field equations for the self-regulating IRP}\label{sec:meanfield}

The usual approach, before actually studying the spatial stochastic model,
is to derive a non-spatial deterministic version called the mean field.
We consider the differential evolution equations for the concentrations
$\{u_i\}_{i \ge 0}$ where $u_i=u_i(t)$ can be thought of as the proportion
of sites with $i$ individuals
at time $t$ when there is a very large population
(clearly $\sum_i u_i=1$). We compute the stationary solutions of the system
of equations: a stationary solution with $u_0<1$ corresponds to
survival, while $u_0 \ge 1$ suggests almost sure extinction.
To avoid unnecessary complications, in this preliminary study, we take $N=1$,
that is, there is no difference between patches and sites.

Let us discuss briefly the mean field equation for the logistic IRP.
The mean field equations are
\[
\begin{cases}
u_0'=u_1-\lambda u_0\sum_{i=1}^\kappa iu_i\\
u_1'=2u_{2}+\lambda u_0\sum_{i=1}^\kappa iu_i-u_1(1+\phi (1-\frac{1}{\kappa})),\\
u_i'=(i+1)u_{i+1}+(i-1)\left(1-\frac{i-1}{\kappa}\right)\phi u_{i-1}-i\left(1+\left(1-\frac i\kappa\right)\phi\right)u_i, & 1\le i\le \kappa-1,\\
u_\kappa'=-\kappa u_\kappa+(\kappa-1)\left(1-\frac{\kappa-1}{\kappa}\right)\phi u_{\kappa-1}.
\end{cases}
\]
Put $u_i'=0$, for $i=1,\ldots,\kappa$ and sum equations from $i$ to $\kappa$.
Then the solutions, for $i=1,\ldots,\kappa$, is
\[
u_i=\frac{(\kappa-1)!}{(\kappa-i)!\, i}\left ( \frac{\phi}{\kappa} \right )^{i-1} u_1.
\]
We plug it into the first equation and require $u_0'=0$ to get
\[\begin{split}
u_0&=\frac{1}{\lambda \sum_{i=0}^{\kappa-1}\frac{(\kappa-1)!}{(\kappa-i-1)!}\left ( \frac{\phi}{\kappa} \right )^{i}}\\ 
& = \frac1\lambda
\left(1+\sum_{i=1}^{\kappa-1}\phi^{i}\prod_{j=1}^i\left(1-\frac{j}{\kappa}\right)\right)^{-1}
\end{split}
\]
From this we get that $u_0\to0$ as $\kappa\to\infty$ (thus indicating the possibility of an endemic
state for $\kappa$ large) when $\phi\ge1$ ($\phi>1$ implies exponential convergence, when $\phi=1$
then $u_0 \sim C/\lambda\sqrt \kappa$).

More precisely, Monotone Convergence Theorem implies that $u_0 \downarrow \min(0, \frac{1-\phi}{\lambda})$
as $\kappa \to \infty$.
Hence, if
$\phi+\lambda \le 1$ then $u_0\ge1$ for all $\kappa$, while if $\phi+\lambda>1$ then $u_0<1$
for $\kappa$ sufficiently large.

If we consider the self-regulating IRP 
the mean field equations are
\[
\begin{cases}
u_0'=u_1-\lambda u_0\sum_{i=1}^\infty iu_i\\
u_1'=2u_{2}+\lambda u_0\sum_{i=1}^\infty iu_i-u_1(1+\phi c(1)),\\
u_i'=(i+1)u_{i+1}+(i-1)c(i-1)\phi u_{i-1}-iu_i(1+\phi c(i)), & 1\le i.
\end{cases}
\]
We can write the solution of the system $u_i'=0$, for $i\ge1$ as
\[
u_i=\frac{\phi^{i-1}}{i}  \prod_{l=0}^{i-1} c(l)     
u_1.
\]
We plug it into the first equation and require $u_0'=0$ to get
\[\begin{split}
u_0&=\frac{1}{\lambda \sum_{i=0}^{\infty}\phi^{i}\prod_{l=0}^i c(l)}.\\
\end{split}
\]
If we are given a sequence of processes regulated by the
functions $\{c_n\}_{n\ge0}$, 
since $c_n \uparrow c_\infty$ as $n \to \infty$, 
by Monotone Convergence Theorem we have
\[
u_0 =\frac{1}{\lambda \sum_{i=0}^{\infty}\phi^{i}\prod_{l=0}^i c_n(l)}
\downarrow \frac{1}{\lambda \sum_{i=0}^{\infty}\phi^{i}\prod_{l=0}^i c_\infty(l)}
\ge \min \left (0, \frac{1-\phi}{\lambda} \right )
\]
and the equality holds if and only if $c_\infty \equiv 1$.
In this case, as before for the logistic IRP,
$\phi+\lambda \le 1$ implies $u_0\ge1$ for all $c_n$, while $\phi+\lambda>1$ implies $u_0<1$
for $n$ sufficiently large.

Thus
the mean field model for the self-regulating IRP suggests that, in the spatial case, if $\phi+2d\lambda \le 1$
there is extinction for all controlling functions $c_n$
(see Theorem~\ref{th:brwcoupling}) and, if
$c_\infty \equiv 1$,  survival is implied by $\phi+2d\lambda>1$ provided that $n$ is sufficiently large
(we prove a slightly different result, see Theorem~\ref{th:cgen}).

\section{Main results and discussion}\label{sec:results}

In this section we discuss the effect of the parameters
$\lambda$, $\phi$, $c$ and $N$ on the behavior of the PRP.
Part of the arguments are done using monotonicity and
coupling with known processes.
Indeed for any fixed $N$, the PRP
is attractive with respect to $\lambda, \phi, c$. Moreover we may couple
a PRP on $\cE^d_N$ with another one on $\cE^d_1\equiv \Z^d$ in the following
natural way.
Consider two PRPs:
$\{\eta_t\}_{t \ge 0}$  on $\cE^d_N$ (with parameters $\lambda, \phi$ and function $c$) and
$\{\xi_t\}_{t \ge 0}$  on $\Z^d$ (with parameters $\lambda_1, \phi_1$ and function $c_1$). If
$\lambda\ge \lambda_1$, $\phi \ge \phi_1$, $c \ge c_1$ then the projection of $\eta_t$ on
$\Zd$ (namely $\sum_{r \in \K_N} \eta_t(\cdot,r)$) dominates $\xi_t(\cdot)$.
%
%

This kind of coupling is the key of the proof of the following result, which
states that if the breeding parameters $\lambda$ and $\phi$ are sufficiently small
then we have almost sure extinction.
\begin{theorem}\label{th:brwcoupling}
For all functions $c$, if $\phi+2d\lambda \le 1$ then there is a.s.~extinction.
\end{theorem}

By monotonicity it is clear that, if $N$, $c$ and $\phi$
are fixed then there exists a critical $\lambda_{cr}=\lambda_{cr}(\phi,c,N)$
such that given $\lambda> \lambda_{cr}$ the PRP survives with positive probability,
while $\lambda<\lambda_{cr}$ implies almost sure extinction. By stochastic domination
$\lambda_{cr} \le \lambda_{CP}(\cE^d_N)$ thus if $\lambda > \lambda_{CP}(\cE^d_N)$
then then PRP survives with positive probability for all choices of $\phi$ and $c$.

The following theorem describes the dependence of $\lambda_{cr}$ on $\phi$ and $c$.
Since $c(i)$ can be equal to $0$ for some $i$,
we identify $1/0$ with $+\infty$.

\begin{theorem}\label{th:lambdagen}
\begin{enumerate}[(1)]
\item
Given any $\phi$, $c$ such that $\sum_{n=0}^\infty (\phi^{n} \prod_{i=0}^n c(i))^{-1}<+\infty$ then
for every $\lambda \ge 0$ the PRP
 survives with positive probability, hence $\lambda_{cr}(\phi,c,N)=0$.
\item
Given any $\phi$, $c$ such that $\sum_{n=0}^\infty \phi^{n} \prod_{i=0}^n c(i)<+\infty$ then 
$
\lambda_{cr}(\phi,c,N)
\in(0,\lambda_{CP}(\cE^d_N))$.
\end{enumerate}
\end{theorem}

In particular from the previous theorem we deduce that the value of
$\lim_{i\to\infty}c(i)$ (which exists by the assumption of monotonicity that we made)
tells us whether $\lambda_{cr}=0$ or not in almost every case.

\begin{cor}\label{cor:lambdagen}
Let $c(\infty)=\lim_{i\to\infty}c(i)$. If $\phi\cdot c(\infty)>1$ then $\lambda_{cr}(\phi,c,N)=0$.
If $\phi\cdot c(\infty)<1$ then $\lambda_{cr}(\phi,c,N)>0$. As a particular case, for all functions
$c$, $\phi<1$ implies $\lambda_{cr}>0$.
\end{cor}

In the case $\phi \cdot c(\infty)=1$ the behavior depends on the speed of convergence
of $c(n)$ to $c(\infty)$. For instance if $c(n):=\frac{(n+3)^2}{2(n+2)^2}$ and $\phi=2$
then Theorem~\ref{th:lambdagen}(1) applies and $\lambda_{cr}=0$. Unfortunately,
due to the monotonicity of $c$, in the case $\phi \cdot c(\infty)=1$,
Theorem~\ref{th:lambdagen}(2) is useless (since the series is always divergent).

\medskip

Again by monotonicity, if $N$, $c$ and $\lambda$
are fixed then there exists a critical $\phi_{cr}=\phi_{cr}(\lambda,c,N)$
such that given $\phi> \phi_{cr}$ the PRP survives with positive probability,
while $\phi<\phi_{cr}$ implies almost sure extinction. In general it could happen that
$\phi_{cr}=\infty$ (almost sure extinction for all $\phi$ - for instance when $N=1$, $c=\delta_0$ and
$\lambda <\lambda_{CP}(\Z^d)$) or
$\phi_{cr}=0$ (survival for all positive $\phi$ - for instance when $\lambda>\lambda_{CP}(\cE^d_N)$).
The following theorem gives sufficient conditions for $\phi_{cr} \in (0, \infty)$.
\begin{theorem}\label{th:phigen}
\begin{enumerate}[(1)]
\item
If $\lambda>\lambda_{CP}(\cE^d_N)$ then $\phi_{cr}(\lambda,c,N)=0$ for all $c$, $N$ and the PRP survives
with positive probability when $\phi=0$.
\item
For any $\lambda\in(0, 1/2d)$ and $c$ such that $c(1)>0$ then $\phi_{cr}(\lambda,c,N) \in (0, +\infty)$.
\end{enumerate}
\end{theorem}

Actually to prove $\phi_{cr}<\infty$ we only need the hypothesis $c(1)>0$ while
$\lambda \in (0,1/2d)$ is needed to prove $\phi_{cr}>0$. One guesses that
$\phi_{cr}>0$ could be proved under the milder assumption $\lambda\in(0,\lambda_{CP}(\cE^d_N))$;
in order to do so one could adapt the proofs of \cite[Theorem 3]{cf:BL} and  of
\cite[Theorem 1(c)]{cf:Schin2002}. The main difficulty is to prove the analog of \cite[Theorem 1.7]{cf:BG91} for
the CP on $\cE^d_N$. Since $\lambda_{CP}(\cE^d_N) \downarrow 1/2d$, extending the interval for $\lambda$ to
$(0,\lambda_{CP}(\cE^d_N))$ seems to be a minor improvement
at least for large values of $N$.

%
%


\medskip

Now we want to study the dependence of the behavior of the PRP
on the underlying space: note that ``space'' here is
described both by $N$ (the horizontal space) and $c$ (the vertical space).
Hence it is natural to define a sequence of PRPs by means of
 nondecreasing sequences of functions $\{c_n\}_{n \in \N}$ and
of patch dimensions $\{N_n\}_{n \in \N}$. More precisely
the $i$th process is $\N^{\cE^d_{N_n}}$-valued and its transition rates are defined
as in  equation \eqref{eq:genrates} with $c=c_n$ (being $d, \lambda, \phi$ fixed for all $n$).

\begin{theorem}\label{th:cgen}
Consider a sequence of PRPs on $\cE^d_{N_n}$
with parameters $\phi$ and $\lambda$ and control functions $c_n$.
If $\lambda>0$ and one of the following conditions holds
\begin{enumerate}[(1)]
\item
$\phi+2d\lambda>1$ and $N_n\to\infty$;
\item $\phi \, \inf_{i\in \N} c_\infty(i) >1$; 
\end{enumerate}
then there exists $n_{cr}=n_{cr}(\lambda,\phi,\mathbf{c})$ such that for all $n \ge n_{cr}$ there is survival with positive probability.
\end{theorem}

Roughly speaking, Theorem~\ref{th:cgen} states that if the space is sufficiently large then the PRP
survives with positive probability provided that the breeding parameters are not too small.
In particular Theorem~\ref{th:cgen}(1) is a partial converse of Theorem~\ref{th:brwcoupling}.
In the proof we need to mimic a technique introduced in \cite{cf:BZ08} in order to
show the ``convergence'' of the projections on $\Z^d$ of the PRPs to a branching random walk.

\begin{rem}\label{rem:contact}
Theorems~\ref{th:brwcoupling} and \ref{th:cgen} imply, roughly speaking, that when
$N_n \to \infty$, the projection of the $n$-th PRP on $\Zd$
behaves, in the limit $n\to\infty$, as the branching random walk on $\Z^d$
with intra-site breeding rate $\phi$ and inter-site rate $\lambda$,
in the sense that it survives eventually if and only if $\phi +2d\lambda>1$.

Recall that by  \cite[Theorem 5]{cf:BL} for the IRP with parameters $\lambda$ and $\phi$ and with $\kappa$ maximal number of particles per vertex,
there is  a critical value for $\kappa$, say $\kappa_c(\lambda,\phi)$, such that for $\kappa\ge\kappa_c(\lambda,\phi)$ there is
survival and for $\kappa<\kappa_c(\lambda,\phi)$ there is extinction.
It is easy to prove that condition (2) of Theorem~\ref{th:cgen} can be relaxed,
provided we have some knowledge of the function $\kappa_c$. More precisely,
one can prove (using the same coupling arguments of the proof of Theorem~\ref{th:cgen}) that
if, given $\lambda >0$, $\phi \ge 0$ and $\{c_n\}$, there exists $\widetilde \phi>1$
such that $\phi \, c_\infty(\kappa_c(\lambda,\widetilde \phi)-1) > \widetilde \phi$ then
there exists $n_{cr}$ such that for all $n \ge n_{cr}$ there is survival with positive probability.
\end{rem}

The following proposition deals with the asymptotic behavior of the critical
parameters of various processes on $\cE^d_{N_n}$
as $N_n\to\infty$. The results follow from Theorems~\ref{th:brwcoupling}
and \ref{th:cgen}.

\begin{pro}\label{th:Ninfty}
If  $N_n \uparrow \infty$ the following hold.
\begin{enumerate}[(1)]
 \item
$\lambda_{CP}(\cE^d_{N_n})\downarrow 1/2d$.
\item
Consider the biparametric contact process, for all fixed $\phi\ge0$ we have that
\[
\lambda_{BCP}(\cE^d_{N_n},\phi)\downarrow \max\left(\frac{1-\phi}{2d},0\right), 
\]
while for all fixed $\lambda>0$ we have that
\[
\phi_{BCP}(\cE^d_{N_n},\lambda)\downarrow \max\left({1-2d\lambda},0\right).
\]
In the particular case $\phi=\lambda$
(the intra-patch reproduction has the same rate as the inter-patch one), 
the critical parameter converges to $1/(2d+1)$.
\item
For a sequence of PRPs,
\[
\begin{split}
&\lim_{n\to\infty} \lambda_{cr}(\phi,c_n,N_n)=0 \qquad \text{if } \phi>1, \\
&\lim_{n\to\infty} \phi_{cr}(\lambda,c_n,N_n)=0 \qquad \text{if } \lambda>1/2d.
\end{split}
\]
\end{enumerate}
\end{pro}
We note that Proposition~\ref{th:Ninfty}(3) gives the limit of $\lambda_{cr}$ as the size of
the patch goes to infinity.
The following proposition deals with the limit of $\lambda_{cr}$ in
the case $c_\infty\equiv1$ (being $N$ and $\phi$ fixed). It is a consequence
of Corollary~\ref{cor:lambdagen} and Theorem~\ref{th:cgen}(2) and establishes the continuity
of $\lambda_{cr}$ with respect to $c$.
\begin{pro}\label{th:cinfty}
 Given $N\in \N$, $\phi>1$ and a sequence of functions $\{c_n\}_{n \ge 0}$ such that
 $c_\infty\equiv1$, then
\[
\lim_{n\to\infty} \lambda_{cr}(\phi,c_n,N)=\lambda_{cr}(\phi,c_\infty,N).
\]
\end{pro}

We recall that the IRP with cluster size $\kappa$ is a particular case
of the PRP (taking $N=1$ and $c=\ident_{\{0,\ldots,\kappa-1\}}$).
Thus if $\bar \lambda^\kappa(\phi)$ is the critical parameter of this IRP,
from Proposition~\ref{th:cinfty} we get the following result on
the asymptotic behavior of $\bar \lambda^\kappa(\phi)$
as $\kappa$ goes to infinity.

\begin{cor}\label{th:irpk}
If $\phi>1$ then
 $\lim_{\kappa \to\infty} \bar \lambda^\kappa(\phi)=0$.
\end{cor}

We note that, if $\bar \lambda^\infty(\phi)$ is the critical parameter of the IRP
with cluster size equal to infinity (see \cite{cf:Be} for some details
on this process), then as a consequence of Corollary~\ref{cor:lambdagen}, if $\phi>1$ then $\bar \lambda^\infty(\phi)=0$.
Thus the previous corollary is a result of continuity at infinity for $\bar \lambda^\kappa(\phi)$  with respect to the cluster size.

\section{Proofs}\label{sec:proofs}

\begin{proof}[Proof of Theorem~\ref{th:brwcoupling}]
It is enough to note that the total reproduction rate of each particle
is bounded from above by $\phi+2d\lambda$.
Indeed the projection
$\sum_{k \in \K_n} \tilde \eta_t(x,k)$
  is dominated by a branching random walk on $\Zd$ with intra-patch infection rate $\phi$
and inter-patch rate $\lambda$ on each edge.
In this branching random walk 
each particle breeds at rate $\phi +2d\, \lambda$ and dies at rate 1.
\end{proof}

In order to prove Theorem~\ref{th:lambdagen} we need the following lemma on discrete time
random walks.
\begin{lem}\label{lem:rw}
Let $\{Z_n\}_{n\ge0}$ be a random walk on $\N^N$ with the following transition probabilities:
\begin{equation}\label{eq:discrete}
 \begin{cases}
(i_1, \ldots, i_N) \longrightarrow (i_1, \ldots, i_j+1, \ldots, i_N) \quad \text{with probability } \frac{\phi c(i_j)}{(1+\phi)N}\\
(i_1, \ldots, i_N) \longrightarrow (i_1, \ldots, i_j-1, \ldots, i_N) \quad \text{with probability } 
\frac{i_j}{(1+\phi) \sum_{r=1}^N i_r}\\
(i_1, \ldots, i_N) \longrightarrow (i_1, \ldots, i_N) \quad \text{with probability } \frac{\phi}{(1+\phi)N} \sum_{j=1}^N (1-c(i_j)).\\
\end{cases}
\end{equation}
\begin{enumerate}[(1)]
\item
If $N=1$ then the random walk is transient if and only if 
$\sum_{n=0}^\infty (\phi^{n} \prod_{i=0}^n c(i))^{-1}<+\infty$.
\item
For any $N$, if $\sum_{n=0}^\infty \phi^{n} \prod_{i=0}^n
c(i)<+\infty$ then the random walk is positive recurrent.
\end{enumerate}
\end{lem}

\begin{proof}
\begin{enumerate}[$(1)$]
\item
If $N=1$ then the random walk is a birth-death process
with forward probabilities $\phi c(i)/((1+\phi))$ and backward probabilities $1/((1+\phi))$
which by \cite[Theorem 5.9]{cf:Woess} is transient if and only if
$\sum_{n=0}^\infty (\phi^{n} \prod_{i=0}^n c(i))^{-1}<+\infty$.
\item
Positive recurrence is equivalent to the existence of a finite invariant measure
which, in this case, is unique up to multiplication.
It is not difficult to prove that
\[
\nu
(i_1, \ldots, i_N)
=
\begin{cases}
\frac{(\sum_{j=1}^N i_j -1)!}{\prod_{j=1}^N i_j!}
 \prod_{j : i_j>0} \left ( \left(\frac{\phi}{N} \right )^{i_j} 
\prod_{i=0}^{i_j-1} c(i)\right), & \text{if } \sum_{j=1}^N i_j >0 \\
\qquad 1, & \text{if } i_j=0, \ \forall j=1, \ldots N
\end{cases}
\]
is a reversible measure. The claim follows noting that, if $\sum_{j=1}^N i_j >0$,
\[
 \nu(i_1, \ldots, i_N) \leq 
\frac{(\sum_{j=1}^N i_j)!}{\prod_{j=1}^N i_j!}\frac{1}{N^{\sum_{j=1}^N i_j}}
 \prod_{j : i_j>0} \left ( \phi^{i_j} 
\prod_{i=0}^{i_j-1} c(i)\right) \le \prod_{j : i_j>0} \left ( \phi^{i_j} 
\prod_{i=0}^{i_j-1} c(i)\right)
\]
whence
\[
\nu(\N^N) = \sum_{i_1,\ldots,i_N}\nu
(i_1, \ldots, i_N)
\le \left(1+\phi\sum_{n=0}^{+\infty}\phi^n\prod_{i=0}^nc(i)
\right)^N < +\infty.
\]
\end{enumerate}
\end{proof}

\begin{proof}[Proof of Theorem~\ref{th:lambdagen}]
\begin{enumerate}[(1)]
\item
Take the PRP with $\lambda=0$. Since it cannot leave the original patch,
this is a (continuous time) random walk on $\N^N$ and it is stochastically dominated by the
original PRP; in this case survival means ``avoiding the origin''.
Consider a new process in the patch which is as the original one but
the success of a reproduction trial in the $j$-th vertex of the patch is $c(\sum_{r=1}^N i_r)$
instead of $c(i_j)$. Thus, the total number of particles in the patch
behaves as a PRP with $N=1$, hence we are left to prove the statement
when $N=1$.
The survival of the (continuous time) random walk is equivalent
to the survival of its discrete counterpart (see equation \eqref{eq:discrete}) which
is guaranteed by Lemma~\ref{lem:rw}(1) noting that transience implies survival.
\item
By (1), $\lambda_{cr} \le \lambda_{CP}(\cE^d_N)$. We now prove that $\lambda_{cr}>0$.
Let us start by considering the case $N=1$ (remember that $\cE^d_1$ can be identified with $\Z^d$).
For $\lambda$ small, we want to find a subcritical branching process which dominates the total number of successful inter-patch reproductions.
Indeed in this case the PRP cannot survive without leaving the first patch, hence survival
is equivalent to the positive probability of having an infinite number of successful
inter-patch infections.
  Once a particle appears at $x$, the total number of particles living at $x$ up to the next time
when the colony at $x$ has no individuals is, according to Lemma~\ref{lem:rw},
a positive recurrent random walk (let us call $\tau_0$ the number of steps before reaching 0) with the following rates:
\[
\begin{cases}
i \longrightarrow  i+1 \quad \text{with rate } i\phi c(i)\\
i \longrightarrow  i-1 \quad \text{with rate } i\\
i \longrightarrow  i \quad \text{with rate } i \phi (1-c(i)).\\
\end{cases}
\]
Note that the exponential clock of the transition (including the loop) has parameter $i(1+\phi)$.

The number of inter-patch reproduction trials made from the site $x$ depends on the number
of particles at $x$: indeed, when the number of particles at $x$ is $i$,
the total reproduction rate towards the neighboring patches is $i \lambda$. Hence
the number of inter-patch reproduction trials between two transitions does not depend on $i$, since it is the number of times a Poisson clock with rate $i \lambda$ rings before the transition clock
with rate $i(1+\phi)$ does. It is clear that
the total number of reproduction trials (before the progeny of the original particle at $x$ dies out) from $x$ to the neighbors is $\sum_{k=1}^{\tau_0} \sum_{j=1}^{2d}Y_k(j)$ where $\{Y_k\}_{k \ge 1}$ is an i.i.d.~sequence
of random $2d$-dimensional vectors; to be precise, if $n_i$ is the number of reproduction trials on
the neighbor $i$ then
\[
\pr(Y_k=(n_1, \ldots,n_{2d}))=\frac{(\sum_{i=1}^{2d} n_i)! (\lambda/(1+\phi))^{\sum_{i=1}^{2d} n_i}}{(1+2d\lambda/(1+\phi))^{1+\sum_{i=1}^{2d} n_i}
\prod_{i=1}^{2d}n_i! }
\]
(see \cite{cf:BZ08bis} for details).
Thus
\[
\E \left ( \sum_{k=1}^{\tau_0} \sum_{j=1}^{2d}Y_k(j) \right )=2d\frac{\lambda}{1+\phi} \E(\tau_0) \to 0
\]
as $\lambda \to 0$. By choosing $\lambda < (1+\phi)/(2d\E(\tau_0))$ the total number of successful trials
(which cannot exceed the total number of trials) is
dominated by a (subcritical) branching process with expected number of offspring $2d\frac{\lambda}{1+\phi} \E(\tau_0)<1$.


Let us take $N>1$: we cannot repeat exactly the same argument since generations inside a patch
are not independent, due to the presence of the controlling function $c$.
Hence we couple with a suitable process. The coupled process has the following rules:
\begin{enumerate}[$a.$]
 \item
 every inter-patch reproduction is successful (i.e.~we do not use $c$);
\item
 we define the particle born from an inter-patch reproduction as the ancestor of
its descendance inside its patch;
\item
 inside a patch particles compete for resources if and only if they have the same ancestor
 (i.e.~the function $c$ applies inside
the same progeny).
\end{enumerate}
Clearly this process dominates the original PRP and progenies of different ancestors
are independent.

The number of generations of the progeny of a particle in the dominating process 
is equal to the absorption time $\tau_0$ of a randow walk with rates
\[
\begin{cases}
(i_1, \ldots, i_N) \longrightarrow (i_1, \ldots, i_j+1, \ldots, i_N) \quad \text{with rate } \frac{\phi}{N} c(i_j)\sum_{r=1}^N i_r\\
(i_1, \ldots, i_N) \longrightarrow (i_1, \ldots, i_j-1, \ldots, i_N) \quad \text{with rate } i_j\\
(i_1, \ldots, i_N) \longrightarrow (i_1, \ldots, i_N) \quad \text{with rate } \frac{\phi}{N}\sum_{j=1}^N (1-c(i_j)) \sum_{r=1}^N i_r.\\
\end{cases}
\]
As before, the behavior of this random walk is equivalent to the behavior of its discrete-time counterpart (see equation \eqref{eq:discrete})
which,
according to Lemma~\ref{lem:rw}, is positive recurrent (thus $\E(\tau_0)<+\infty$).
Hence, being the discrete time counterpart recurrent,
 survival is equivalent to the positive probability of having an infinite number of progenies.
But 
the expected number of successful inter-patch reproductions is equal to
 $2d\frac{\lambda}{1+\phi} \E(\tau_0)$. Choosing $\lambda$ small enough
yields the conclusion.
\end{enumerate}
\end{proof}

\begin{proof}[Proof of Theorem~\ref{th:phigen}]
\begin{enumerate}[(1)]
\item
It follows from the fact that for all $\phi \ge 0$ the PRP stochastically
dominates the CP on $\cE^d_N$.
\item
The hypothesis $c(1)>0$ yields $\phi_{cr}<\infty$. Indeed
let $\tilde \phi_{cr}(\lambda,2)>0$ be the critical threshold in \cite[Theorem 3]{cf:BL} and
choose $\phi > \tilde \phi_{cr}(\lambda,2)/c(1)$.
Let $\widetilde\eta_t:\Z^d\to\{0,1,2\}$ be the IRP with parameters $\widetilde\lambda=\lambda$,
$\widetilde\phi=\phi c(1)$ and $\widetilde\kappa=2$. By these choices of the parameters,
$\widetilde\eta_t$ is supercritical. On the other hand,
the projection  $\xi_t(x)=\sum_{i \in \K_n} \tilde \eta_t(x,i)$
dominates $\widetilde\eta_t$ and this proves that the PRP survives
when $\phi > \tilde \phi_{cr}(\lambda,2)/c(1)$.

To show that $\phi_{cr}>0$ in the case where $\lambda<1/2d$, note that
by Theorem~\ref{th:brwcoupling} if we choose a positive $\phi$ smaller than $1-2d\lambda$
then the population dies out.
\end{enumerate}
\end{proof}

%

\begin{proof}[Proof of Theorem~\ref{th:cgen}]
\begin{enumerate}[(1)]
\item
The PRPs of the sequence dominate the ones where the control functions are equal to
$\delta_0$ for all $n$. Hence it is enough to prove the statement for
the latter case.

The projection of this process on $\Z^d$, namely $\xi^n_t(x):=\sum_{r \in \K_N}\eta^n_t(x,r)$, has transition rates at $x \in \Z^d$
\begin{equation}
\begin{split}
j\to j+1& \text{ at rate }\left (j\phi+\lambda \sum_{y \sim x}\xi^n_t(y) \right )
\left (1-\frac{j}{N_n} \right );\\
j\to j-1 & \text{ at rate }j.
\end{split}
\end{equation}
%
%
This
process dominates eventually (as $N_n \to \infty$) the process $\{\widetilde \xi_t\}_{t \ge 0}$ with rates
\begin{equation}\label{eq:genrates2}
\begin{split}
j\to j+1& \text{ at rate }\left (j\phi+\lambda \sum_{y \sim x} \widetilde \xi_t(y) \right )(1-\eps)\ident_{[0,\bar N]}(j);\\
j\to j-1 & \text{ at rate }j;
\end{split}
\end{equation}
for all $\bar N>0$ and $\eps >0$ (it suffices that $N_n>\bar N/\eps$).
Choose $\eps>0$
such that $(\phi+2d \lambda)(1-\eps)>1$. We will prove now, following \cite{cf:BZ08},
that for all sufficiently large $\bar N$, the process   $\{\widetilde \xi_t\}_{t \ge 0}$ survives with
positive probability.
The strategy is to study the branching random walk
$\{\bar \xi_t\}_{t \ge 0}$ with rates
\begin{equation}\label{eq:genrates3}
\begin{split}
j\to j+1& \text{ at rate }\left (j\bar \phi+\bar\lambda \sum_{y \sim x}\bar \xi_t(y) \right ),\\
j\to j-1 & \text{ at rate }j,
\end{split}
\end{equation}
starting with one particle at the origin $0 \in \Z^d$, where $\bar \phi=(1-\eps)\phi$ and $\bar \lambda=(1-\eps)\lambda$. Note that the process $\{\widetilde \xi_t\}_{t \ge 0}$ can be seen as a truncation
(at height $\bar N$) of $\{\bar \xi_t\}_{t \ge 0}$.

It is not difficult to verify that $\E^{\delta_0}(\bar \xi_t)$ satisfies the
following differential equation (see \cite[Section 5]{cf:BZ08}, where we did the same with
the classical BRW):
\begin{equation}
 \frac{d}{dt}\E^{\delta_0}(\bar \xi_t(x))=-\E^{\delta_0}(\bar \xi_t(x))+\phi\E^{\delta_0}(\bar \xi_t(x))
+\lambda\sum_{\stackrel{y\sim x}{y\neq x}}\E^{\delta_0}(\bar \xi_t(y)),
\end{equation}
whose solution is
\begin{equation}\label{eq:noimm}
\E^{\delta_0}(\bar \xi_t(x))=
\sum_{n=0}^\infty \sum_{k=0}^{n-|x|_{\Z^d}} \mu^{(n,k)}(0,x)
\frac{{\bar \phi}^{k}{\bar \lambda}^{n-k} t^n}{n!}e^{-t},
\end{equation}
where $\mu^{(n,k)}(0,x)$ is the number of paths from $0$ to $x$ of length $n$ and $k$ loops.
Moreover, if $|x|_{\Z^d}=1$ then, taking $t=n$ large enough,
\[
\begin{split}
\E^{\delta_0}(\bar \xi_t(x)) &\ge
\sum_{k=0}^{n-1} \mu^{(n,k)}(0,x)\frac{{\bar \phi}^{k}{\bar \lambda}^{n-k} t^n}{n!}e^{-t}=
\sum_{k=0}^{n-1} \mu^{(n,k)}(0,x)\frac{{\bar \phi}^{k}{\bar \lambda}^{n-k}}{(\bar \phi+2d\bar \lambda)^n} \frac{n^n (\bar \phi+2d\bar \lambda)^n}{n!}e^{-n} \\
& \, \stackrel{n \to \infty}{\sim}\, \frac{(\bar \phi+2d\bar \lambda)^n}{\sqrt{2 \pi n}} \sum_{k=0}^{n-1} \mu^{(n,k)}(0,x)\frac{{\bar \phi}^{k}{\bar \lambda}^{n-k}}{(\bar \phi+2d\bar \lambda)^n}
\ge \frac{(\bar \phi+2d\bar \lambda)^n}{\sqrt{2 \pi n}} C n^{-d/2}
\end{split}
\]
for some $C=C(\bar\lambda,\bar\phi)>0$, since $\sum_{k=0}^{n-1} \mu^{(n,k)}(0,x){\bar \phi}^{k}{\bar \lambda}^{n-k}/(\bar \phi+2d\bar \lambda)^n$
is the probability of being in $x$ ($|x|_{\Z^d}=1$) after $n$ steps for a discrete-time random walk with transition probabilities
\[
p(x,y)=
\begin{cases}\displaystyle
\frac{\bar \lambda}{\bar \phi+2d\bar \lambda} & \text{if } x\sim y, \\
\displaystyle
\frac{\bar \phi}{\bar \phi+2d\bar \lambda} & \text{if } x=y, \\
\end{cases}
\]
(see \cite[Corollary 13.11]{cf:Woess2} for the asymptotic estimates of the $n$-step probabilities).
Hence there exists $t=n$ sufficiently large such that $\E^{\delta_0}(\bar \xi_t(x))>1$ (this result is analogous
to Lemmas~5.2 and 5.3 of \cite{cf:BZ08}). Using the same arguments as in Lemma~5.4, Remark~5.5 and Theorem~5.6 of
\cite{cf:BZ08},
it is possible to prove that the process $\{\widetilde \xi_t\}_{t \ge 0}$ survives when $\bar N$ is sufficiently large and this implies
the survival of  $\{\xi^i_t\}_{t \ge 0}$ when $i$ is sufficiently large.
\item

Let $\kappa_c=\kappa_c(\lambda,(1+\phi)/2)$ be the critical value given by  \cite[Theorem 5]{cf:BL}
(see Remark~\ref{rem:contact} for details).
Let $\alpha:=\inf_{i\in \N}c_\infty(i)$: by hypothesis $\alpha>0$ and $\phi\alpha>1$. 
Choose $\varepsilon>0$ such that
$\widetilde\phi:=\phi(\alpha-\varepsilon)>1$. By the definition of $c_\infty$ and the monotonicity 
of the sequence
of functions $c_n$, there exists $n_c$ such that $c_n(\kappa_c(\lambda,\widetilde\phi)-1)>\alpha-\varepsilon$ for all
$n\ge n_c$. Thus for all $n\ge n_c$, if in a patch there are at most $\kappa_c(\lambda,\widetilde\phi)$
 particles, then for the PRP the intra-patch reproduction rate is at least
$\widetilde\phi$. Hence this PRP  dominates the IRP with parameters $\lambda$,
 $\widetilde\phi$ and $\kappa_c(\lambda,\widetilde\phi)$
 which survives with positive probability.
\end{enumerate}
\end{proof}

\section{Open questions}\label{sec:open}

In the phase diagram the answer to some questions remains unknown.
For instance, we do not know if $\lambda_{cr}(\phi,c,N)$ is always equal to zero
when $\phi \cdot c(\infty)=1$ (we know there are examples
where this happens). Moreover we do not know
if $\phi_{cr}(\lambda,c,N)>0$ when $\lambda \in (1/2d, \lambda_{CP}(\cE^d_N))$.
One interesting question is whether a nontrivial invariant
measure exists when there is survival with positive probability.
The answer is positive if $c(i)=0$ eventually (in this case the
configuration space is compact and the invariant measure is obtained
as the limit of the process starting from
the maximal configuration). In the noncompact case one could try
to mimic the technique employed in \cite{cf:BPZ07}.

On the other hand, one
may consider other models which, although similar to those
considered here, cannot be obtained
as particular cases of the PRP.
One of particular interest is a self-regulating IRP where the inter-cluster
infection is always active.
In each site of $\Z^d$ there is a (possibly infinite) cluster and,
given a regulating function $c$, the transition rates at site $x$ are:
\[
\begin{split}
i\to i+1& \text{ at rate }\left(i\phi+\lambda\sum_{z\sim x}\eta(z)\right)c(i),\text{ for }0\le i;\\
i\to i-1 & \text{ at rate }i,\text{ for }i\ge1.
\end{split}
\]
In the present work, we dealt with the case $c(i)=\max(0,1-i/N)$ (see Subsection~\ref{sec:logistic}).
Further investigations may be made on the processe where we
add a rate $\delta$ for the transition $i\to0$
- the so-called \textit{catastrophe} or \textit{mass extinction}.

%
%
%

\section*{Acknowledgments}
The authors acknowledge support from the Italian Research Project
\textit{Prin 2006 - Modelli stocastici a molti gradi di libert\`a: teoria e applicazioni}.
The first author is grateful to the \textit{Dipartimento di Matematica e Applicazioni - Universit\`a
di Milano-Bicocca} for the hospitality and logistic support.
Moreover, the authors want to thank the anonymous referee for
many useful suggestions that helped improving the paper.

\end{document}